\documentclass[12pt]{amsart}
\usepackage{hyperref}
\usepackage{amscd}
\usepackage{mathbbold, color}
\usepackage{mathbbold, cite}
\usepackage{epsfig}
\usepackage{mathdots}
\usepackage{amssymb}
\usepackage{amsfonts}
\usepackage{amsbsy}
\usepackage{graphicx}
\usepackage{enumitem}
\usepackage{ulem}
\newtheorem{theo}{Theorem}[section]
\newtheorem{lem}[theo]{Lemma}
\newtheorem{prop}[theo]{Proposition}

\newtheorem{que}{Question}

\theoremstyle{definition}

\newtheorem{exam}[theo]{Example}
\newtheorem{rem}[theo]{Remark}

\numberwithin{equation}{section}

\newcommand{\abs}[1]{\lvert#1\rvert}

\newcommand{\R}{{\mathbb R}}
\newcommand{\Z}{{\mathbb Z}}

\newcommand{\N}{{\mathbb N}}

\newcommand{\D}{{\mathcal D}}

\newcommand{\eproof}{\hfill$\square$}

\begin{document}
\baselineskip 16pt

\title[]{Spectral measures with arbitrary dimensions}


\author[Yu-Liang Wu]{Yu-Liang Wu}
\address[Yu-Liang Wu]{Department of Mathematical Sciences\\
	University of Oulu, P.O. Box 3000, 90014, Oulu, Finland}
\email{yu-liang.wu@oulu.fi}
\author[Zhi-Yi Wu]{Zhi-Yi Wu}
\address[Zhi-Yi Wu]{Department of Mathematical Sciences\\
	University of Oulu, P.O. Box 3000, 90014, Oulu, Finland}
\email{zhiyiwu@126.com}

%
\date{}



\begin{abstract}
It is known [Dai and Sun, J. Funct. Anal. 268 (2015), 2464--2477] that there exist spectral measures with arbitrary Hausdorff dimensions, and it is natural to pose the question of whether similar phenomena occur for other dimensions of spectral measures. In this paper, we first obtain the formulae of Assouad dimension and of lower dimension for a class of Moran measures in dimension one that is introduced by An and He [J. Funct. Anal. 266 (2014), 343--354]. Based on these results, we show the existence of spectral measures with arbitrary Assound dimensions $\dim_A$ and lower dimensions $\dim_L$ ranging from $0$ to $1$, including non-atomic zero-dimensional spectral measures and one-dimensional singular spectral measures, and prove that the two values may coincide. In fact, more is obtained that for any $0 \leq t \leq s \leq r \leq u\leq 1$, there exists a spectral measure $\mu$ such that
\[\dim_L \mu=t, \dim_H \mu=s, \dim_P\mu=r~ \text{and} \dim_A\mu=u,\]
where $\dim_H$ and $\dim_P$ denote the Hausdorff dimension and packing dimension of the measure $\mu$, respectively. This result improves and generalizes the result of Dai and Sun more simply and flexibly.

\end{abstract}

\maketitle

\section{Introduction}
In classical Fourier analysis, the collection of exponential functions $\{e^{-2\pi inx}:n\in\Z\}$ is known to form an orthonormal basis for the $L^2$ space of Lebesgue measure restricted on $[0,1]$. This idea then gives rise to what defines a \emph{spectral measure}; namely, a compactly supported probability measure on $\R^d$ whose $L^2$ space admits a countable orthonormal basis of the form $\big\{e^{-2\pi i \langle \lambda, x \rangle}:\lambda\in\Lambda\big\}$ with $\Lambda$ a subset of $\R^d$. In 1998, Jorgensen and Pederson \cite{JP}, found that the standard fourth-middle Cantor measure is a spectral measure while the usual middle-third Cantor measure is not. Thereafter, there has been a deal of activity concerning the study of spectrality of fractal measures, including self-similar measures, self-affine measures, and Moran measures (cf. \cite{LW02,AFLai,DHLai19,DHLai13,Dai16}). In the meanwhile, discovered are many exotic phenomena of singular spectral measures which are not observed among classical measures \cite{DS15,DHLai13,FHW18,LW02}. Among these, Dai and Sun \cite{DS15} proved the following somewhat surprising result for one-dimensional singular spectral measures. 
\begin{theo} \label{thm:hausdorff_intermediate}
	For any $t\in[0,1]$, there exists a singular spectral measure whose Hausdorff dimension is $t$.
\end{theo}
\noindent This theorem then inspires the following question:
\begin{que} \label{question_intermediate}
    Is the intermediate property above seen for other varieties of dimension?
\end{que}
\noindent In this paper, we show that, in addition to Hausdorff dimension, the answer to the question is affirmative for packing dimension, Assouad dimension, and lower dimension. It is known (see for example \cite{F21}) that any compactly supported Borel probability measure $\mu$ satisfies 
\begin{align}\label{eqdim}
    \dim_L\mu \le \dim_H\mu \le \dim_P\mu \le \dim_A\mu,
\end{align}
where $\dim_L$, $\dim_H$, $\dim_P$ and $\dim_A$ denote the lower, Hausdorff, packing and Assouad dimension of $\mu$, respectively. This fact then motivates the main theorem of the paper which addresses Question \ref{question_intermediate}.
\begin{theo} \label{thm:arbitrary_dimension}
	For any $0\leq t_L\leq t_H\leq t_P\leq t_A\leq 1$, there exists a spectral measure such that 
	\begin{equation}
		\dim_L \mu=t_L, \dim_H \mu=t_H, \dim_P\mu=t_P~ \text{and}~ \dim_A\mu=t_A,
	\end{equation}
\end{theo}
\noindent It is to be noted that not only this result improves Theorem \ref{thm:hausdorff_intermediate}, but our arguments are simpler and more flexible. 

Our construction of a measure satisfying Theorem \ref{thm:arbitrary_dimension} is based on the following class of Moran spectral measures, which are first considered by An and He \cite{AH14}. Let $\mathcal{B}=\{b_n\}_{n=1}^\infty$ and $\mathcal{D}=\{q_n\}_{n=1}^\infty$ be two sequences of integers with $2\le q_n<b_n$. For each $k\ge 1$, write $\mathcal{D}_k=\{0,1,\ldots, q_k-1\}$ and let
\begin{equation}\label{MM}
	\mu=\mu_{\mathcal{D},\mathcal{B}}:=\delta_{b_1^{-1}\mathcal{D}_1}\ast \delta_{(b_1b_2)^{-1}\mathcal{D}_2}\ast\cdots \ast \delta_{(b_1b_2\cdots b_k)^{-1}\mathcal{D}_k}\ast\cdots,
\end{equation}
where $\delta_{\mathcal{E}}=\frac{1}{\#\mathcal{E}}\sum_{e\in \mathcal{E}}\delta_e$, $\#\mathcal{E}$ denotes the cardinality of the set $\mathcal{E}$, $\delta_e$ denotes the Dirac measure at $e$, and the convergence is in the weak sense.
We call such $\mu$ a \emph{Moran measure}, which becomes a \emph{Cantor measure} $\mu_{q, b}$ when $b_n=b$ and $q_n=q$ for all $n \ge 1$ as is studied in \cite{DHLai13}. The support of $\mu$ can be represented in terms of digit expansion
\begin{equation}\label{supp}
	\textrm{supp}(\mu)=\left\{\sum_{k=1}^\infty \frac{d_k}{b_1\ldots b_k}:~d_k\in \mathcal{D}_{k}, k\ge 1\right\}
\end{equation}
and can be regarded as a special Moran set introduced by Feng, Wen and Wu \cite{FWW}. These Moran measures have been widely used in multifractal analysis \cite{FL}.

The constructions of $\mu$ and $\textrm{supp}(\mu)$ above could be interpreted geometrically as assignment of mass to subintervals of the unit interval $[0,1]$. To begin with, we first equally divide $[0,1]$ into $b_1$ closed subintervals with disjoint interiors. Then, according to the digit set $\D_{1}$, we keep the first $q_{1}$ subintervals, which are called the \emph{intervals of rank $1$}, and uniformly distribute the mass to them; namely, if $I_1$ is an interval of rank $1$ then $\mu(I_1)=1/q_1$. Iteratively, we equally divide each interval of rank $1$ into $b_2$ subintervals with disjoint interiors. According to the digit set $\D_{2}$, we keep the first $q_{2}$ subintervals, call them \emph{intervals of rank $2$}, and uniformly assign to each of them the mass $1/(q_1q_2)$. Continuing the above procedure yields a Borel probability measure $\mu$ whose support $T$ is the intersection of the union of all intervals of rank $k$. Moreover, by Kolmogorov's consistency theorem we have
\[
\mu(I_n)=\frac{1}{q_1\ldots q_n}
\]
for any interval of rank $n\ge 1 $, $I_n$.

To prove Theorem \ref{thm:arbitrary_dimension}, we start by introducing the formulae of the Assouad dimension and the lower dimension of $\mu$. The \emph{Assouad dimension} and the \emph{lower dimension} of a Borel measure $\mu$ on $\R^d$ are two dimensions that have attracted much attention and have been widely studied in the field of fractal geometry, which are defined as 
\begin{equation*} \label{def_assouad}
\begin{split}
 \dim_A \mu=\inf\Bigg\{s\geq0&:~ \text{there exists a constant $C>0$ such that}\\
	&\text{for all $0<r<R<|\textrm{supp}(\mu)|$ and $x\in\textrm{supp}(\mu)$,}\\
	&\quad\quad\quad\quad\quad\quad\quad\quad\quad\quad\quad\quad\frac{\mu(B(x,R))}{\mu(B(x,r))}\leq C\bigg(\frac{R}{r}\bigg)^s\Bigg\}
	\end{split}
\end{equation*}
and
\begin{equation*} \label{def_lower}
\begin{split}
 \dim_L \mu=\sup\Bigg\{s\geq0&:~ \text{there exists a constant $C>0$ such that}\\
	&\text{for all $0<r<R<|\textrm{supp}(\mu)|$ and $x\in\textrm{supp}(\mu)$},\\
	&\quad\quad\quad\quad\quad\quad\quad\quad\quad\quad\quad\quad\frac{\mu(B(x,R))}{\mu(B(x,r))}\geq C\bigg(\frac{R}{r}\bigg)^s\Bigg\}.
	\end{split}
\end{equation*}
Roughly speaking, they are measurements of the optimal global control on the relative measure of concentric balls. The reader may refer to \cite{F21} for more backgrounds and to \cite{HMZ20,HT21} for recent research topics. According to these definitions, the Assouad dimension and lower dimension of $\mu$ can be expressed in terms of sequences $\mathcal{B}$ and $\mathcal{D}$ as follows.
\begin{prop}\label{ass}
Let $\mu$ be the Moran measure defined in \eqref{MM}. Then
\begin{align*}
	\dim_A\mu&=\limsup_{q_kb_{k+1}\cdots b_{k+n}\to+\infty}\frac{\log q_kq_{k+1}\cdots q_{k+n}}{\log q_k b_{k+1}\cdots b_{k+n}}\\
	&:=\lim_{N\to+\infty}\sup_{q_kb_{k+1}\cdots b_{k+n}\geq N}\frac{\log q_k q_{k+1}\cdots q_{k+n}}{\log q_k b_{k+1}\cdots b_{k+n}}.
	\end{align*}
Equivalently,
\begin{enumerate}[label=(\roman*), font=\normalfont]
    \item If $\{q_n\}_{n=1}^{\infty}$ is bounded, then
$$\dim_A\mu=\limsup_{n\to+\infty}\sup_{k\geq1}\frac{\log q_kq_{k+1}\cdots q_{k+n}}{\log b_kb_{k+1}\cdots b_{k+n}};$$
    \item If $\{q_n\}_{n=1}^{\infty}$ is unbounded, then $\dim_A\mu=1$.
\end{enumerate}
\end{prop}
\begin{prop}\label{ass1}
	Let $\mu$ be the Moran measure defined in \eqref{MM}. Then, 
	\begin{align*}
	\dim_L\mu&=\liminf_{b_{k+1}\cdots b_{k+n}\frac{1}{q_{k+n}}\to+\infty}\frac{\log q_{k+1}q_{k+2}\cdots q_{k+n}\frac{1}{q_{k+n}}}{\log b_{k+1}b_{k+2}\cdots b_{k+n}\frac{1}{q_{k+n}}}\\
	&:=\lim_{N\to+\infty}\inf_{b_{k+1}\cdots b_{k+n}\frac{1}{q_{k+n}}\geq N}\frac{\log q_{k+1}\cdots q_{k+n}\frac{1}{q_{k+n}}}{\log b_{k+1}\cdots b_{k+n}\frac{1}{q_{k+n}}}.
	\end{align*}
	In particular,
    \begin{enumerate}[label=(\roman*), font=\normalfont]
        \item If $\{b_n\}_{n=1}^{\infty}$ is bounded, then
$$\dim_L\mu=\liminf_{n\to+\infty}\inf_{k\geq1}\frac{\log q_kq_{k+1}\cdots q_{k+n}}{\log b_kb_{k+1}\cdots b_{k+n}};$$
	   \item If $\{b_n\}_{n=1}^{\infty}$ is unbounded and $\liminf_{n\to+\infty} \frac{q_n}{b_n^{1-\epsilon}}=0$ for any $\epsilon\in[0,1)$, then $\dim_L\mu=0$.
    \end{enumerate}
\end{prop}

\begin{rem}
The limit superior and limit inferior notations used in Propositions \ref{ass} and \ref{ass1} are actually limit superior and limit inferior of a net or Moore–Smith sequence, which is a generalization of sequence first introduced in \cite{Moore1922}. 
\end{rem}

\noindent Given that An and He has proved in \cite{AH14} that $\mu$ is a spectral measure if $q_n$ divides $b_n$ for all $n\ge 1$, we are enabled to show the following weaker version of Theorem \ref{thm:arbitrary_dimension}.
\begin{theo}\label{the-main}
    For any $t\in[0,1]$, there exists a Moran spectral measure $\mu$ such that
    \[\dim_L \mu=\dim_H \mu=\dim_P\mu=\dim_A\mu=t.
    \]
\end{theo}
\noindent It is based on this theorem that the proof of Theorem \ref{thm:arbitrary_dimension} is carried out.
\begin{exam} \label{example}
    Let $q_n=2$  for $n\geq1$. Let $b_n=4^n$ for $n\geq1$. Then $q_n|b_n$ for $n\geq1$ and thus $\mu_{\D,\mathcal{B}}$ is a spectral measure. It is easy to see that the assumption of Proposition \ref{ass1}(ii) is satisfied and thus we have $\dim_L \mu=0$.
\end{exam}

This paper is organized as follows. Propositions \ref{ass} and \ref{ass1} are proved in Section 2 while Theorems \ref{the-main} and \ref{thm:arbitrary_dimension} are proved in Section 3.

\section{Proofs of Propositions \ref{ass} and \ref{ass1}}
In this section, we shall give the proofs of  Propositions \ref{ass} and \ref{ass1}.

\noindent{\bf Proof of Proposition \ref{ass}.} For convenience of discussion, we denote $b_0=q_0=1$ and $s=\limsup_{q_kb_{k+1}\cdots b_{k+n}\to+\infty}\frac{\log q_kq_{k+1}\cdots q_{k+n}}{\log q_k b_{k+1}\cdots b_{k+n}}$. We first prove $\dim_A\mu \geq s$; namely, $\dim_A\mu \geq t$ for every $t < s$. Consider sequences of positive integers $\{k_l\}_{l=1}^{\infty}$ and $\{n_l\}_{l=1}^{\infty}$ such that 
$$\lim_{l \to +\infty} \frac{\log q_{k_l} q_{k_l+1} \cdots q_{k_l+n_l}}{\log q_{k_l} b_{k_l+1} \cdots b_{k_l+n_l}} = s \text{ and } \lim_{l \to +\infty} q_{k_l} b_{k_l+1} b_{k_l+2} \ldots b_{k_l+n_l} =+\infty.$$
Let $q_{k_l} I_{k_l}=\left[0,\frac{q_{k_l}}{b_0 b_1 \cdots b_{k_l}}\right]$ and $I_{k_l+n_l}=\left[0,\frac{1}{b_0 b_1 \cdots b_{k_l+n_l}}\right]$. Then, by the construction of $\mu$, we have that $\mu(q_{k_l}I_{k_l})=\frac{q_{k_l}}{q_1\cdots q_{k_l}}$, that $\mu(I_{k_l+n_l})=\frac{1}{q_1\cdots q_{k_l+n_l}}$, and that
$$\frac{\mu(q_{k_l} I_{k_l})}{\mu(I_{k_l+n_l})}\bigg/\left(\frac{|q_{k_l} I_{k_l}|}{|I_{k_l+n_l}|}\right)^t=\left(q_{k_l}b_{k_l} \cdots b_{k_l+n_l} \right)^{\frac{\log q_{k_l} q_{k_l+1} \cdots q_{k_l+n_l}}{\log q_{k_l} b_{k_l+1} \cdots b_{k_l+n_l}}-t}\to + \infty$$
as $l\to +\infty$. Therefore, by definition of Assouad dimension, we have $\dim_A\mu\geq t$.

On the other hand, we prove $\dim_A\mu \leq s$. Consider any $0<r<R<|\text{supp} \mu|$ and $x\in \textrm{supp}(\mu)$. Choose 
integers $k, l, u \in \N$ and $n \ge 0$ such that 
$$\frac{l}{b_0 b_1 \cdots b_k} < R \leq \frac{l+1}{b_0 b_1 \cdots b_{k}} \leq \frac{2 l}{b_0 b_1 \cdots b_{k}}, \quad 1 \le l < b_k$$
and
$$\frac{u}{b_0 b_1 \cdots b_{k+n}} < r \leq \frac{ u+1}{b_0 b_1 \cdots b_{k+n}} \leq \frac{2 u}{b_0 b_1 \cdots b_{k+n}}, \quad 1 \le u < b_{k+n}.$$ 
Consequently, 
$$\frac{\min\{l,q_k\}}{q_0 q_1 \cdots q_k}\leq \mu(B(x,R))\leq \frac{2 \min\{l+1,q_k\}}{q_0 q_1 \cdots q_k} \le \frac{4 \min\{l,q_k\}}{q_0 q_1 \cdots q_k}$$
and
$$\frac{\min\{u,q_{k+n}\}}{q_0 q_1 \cdots q_{k+n}} \leq \mu(B(x,r))\leq \frac{2 \min\{ u+1,q_{k+n}\}}{q_0 q_1 \cdots q_{k+n}} \leq \frac{4 \min\{u,q_{k+n}\}}{q_0 q_1 \cdots q_{k+n}}.$$
This yields that for any $t>s$, 
\begin{equation}\label{eqleqa}
\begin{aligned}
	& \frac{\mu(B(x,R))}{\mu(B(x,r))}\bigg/\bigg(\frac{R}{r}\bigg)^t \leq \frac{\frac{4 \min\{l,q_k\}}{\min\{u,q_{k+n}\}} \cdot q_{k+1} q_{k+2} \cdots q_{k+n}}{(\frac{l}{2 u} \cdot b_{k+1} b_{k+2} \cdots b_{k+n})^t} \\
	\le& \begin{cases}
	8 \cdot \frac{q_{k} q_{k+1} q_{k+2} \ldots q_{k+n}}{\left(q_{k} b_{k+1} b_{k+2} \ldots b_{k+n}\right)^t} & \text{if } t \le 1 \text{ and } 1 \le u < q_{k+n}; \\
	8 \cdot \frac{q_{k} q_{k+1} \ldots q_{k+n-1}}{\left(q_{k} b_{k+1} \ldots b_{k+n-1}\right)^t} & \text{if } t \le 1 \text{ and } q_{k+n} \le u < b_{k+n}; \\
	8 & \text{if } t > 1.
	\end{cases}
\end{aligned}
\end{equation}
For $1 \ge t>s$, there exists by definition an integer $N=N(t)$ such that
$$\frac{\log q_{k} q_{k+1} \cdots q_{k+l}}{\log q_{k} b_{k+1} \cdots b_{k+l}} < t \text{ whenever } {q_k b_{k+1} \cdots b_{k+n} } > N.$$
Moreover,  if $q_kb_{k+1}\cdots b_{k+n}\leq N$, $\frac{\log q_{k} q_{k+1} \cdots q_{k+l}}{\log q_{k} b_{k+1} \cdots b_{k+l}}<1$. 
Hence, equation \eqref{eqleqa}  is no larger than $\max\{8,8\frac{N}{N^t}\}\leq 8N$. This proves $\dim_A\mu\leq t$, which implies $\dim_A\mu\leq s$ since $t > s$ is arbitrary.

It is left to show that \begin{equation} \label{eq:alt_expression}
	\lim_{N\to+\infty} \sup_{n \ge N, k\geq1}\frac{\log q_k q_{k+1}\cdots q_{k+n}}{\log b_k b_{k+1}\cdots b_{k+n}} = \lim_{N \to +\infty} \sup_{q_k b_{k+1}\cdots b_{k+n} \ge N}\frac{\log q_{k}q_{k+1}\cdots q_{k+n}}{\log q_{k}b_{k+1}\cdots b_{k+n}},
\end{equation}when $\{q_n\}$ is bounded. Since the right-hand side dominates the the left-hand side by definition, we need only to show the remaining inequality.
Now note that if $k_m$ and $n_m$ are sequences such that $q_{k_m} b_{k_m+1}\cdots b_{k_m+n_m} \to \infty$ and that
\[
\lim_{m\to +\infty} \frac{\log q_{k_m} q_{k_m+1}\cdots q_{k_m+n_m}}{\log q_{k_m} b_{k_m+1}\cdots b_{k_m+n_m}}=\lim_{N \to +\infty} \sup_{q_k b_{k+1}\cdots b_{k+n} \ge N}\frac{\log q_{k}q_{k+1}\cdots q_{k+n}}{\log q_{k}b_{k+1}\cdots b_{k+n}},
\]
then $n_m$ is bounded if and only if the limit of the sequence is $0$. Hence, the left-hand side of \eqref{eq:alt_expression} is no less than the right-hand side. 

If $\{q_n\}_{n=1}^{\infty}$ is unbounded, then there exists a subsequence $\{q_{n_k}\}_{k=1}^{\infty}$ of $\{q_n\}_{n=1}^{\infty}$ such that $q_{n_k} \to +\infty$ as $k\to+\infty$. For any $n\geq1$, let $I_n=\Big[0,\frac{1}{b_0 b_1 \cdots b_n}\Big]$ and $q_nI_n=\Big[0,\frac{q_n}{b_0 b_1 \cdots b_n}\Big]$. Then, 
$$\frac{\mu(q_{n_k}I_{n_k})}{\mu(I_{n_k})}=q_{n_k} \text{ and } \left(\frac{|q_{n_k}I_{n_k}|}{|I_{n_k}|}\right)^t=(q_{n_k})^t \quad \text{for all } t \geq 0.$$
Therefore, when $0\leq t< 1$, as $k\to +\infty$
$$\frac{\mu(q_{n_k}I_{n_k})}{\mu(I_{n_k})}\bigg/\left(\frac{|q_{n_k}I_{n_k}|}{|I_{n_k}|}\right)^t=(q_{n_k})^{1-t} \to+\infty.$$
Hence, $\dim_A\mu=1$. 
\eproof

The proof of Proposition \ref{ass1} is similar to that of Proposition \ref{ass}. Nevertheless, we provide the proof for completeness.

\noindent{\bf Proof of Proposition \ref{ass1}.} For convenience of discussion,  we denote $b_0=q_0=1$ and   $s=\lim_{N\to+\infty}\inf_{b_{k+1}\cdots b_{k+n}\frac{1}{q_{k+n}}\geq N}\frac{\log q_{k+1}\cdots q_{k+n}\frac{1}{q_{k+n}}}{\log b_{k+1}\cdots b_{k+n}\frac{1}{q_{k+n}}}$. We first prove $\dim_L\mu \leq s$; namely, $\dim_L\mu \leq t$ for every $t > s$. Consider sequences of positive integers $\{k_l\}_{l=1}^{\infty}$ and $\{n_l\}_{l=1}^{\infty}$ such that 
$$\lim_{l \to +\infty} \frac{\log  q_{k_l+1} \cdots q_{k_l+n_l} \frac{1}{q_{k_l+n_l}}}{\log  b_{k_l+1} \cdots b_{k_l+n_l} \frac{1}{q_{k_l+n_l}}} = s \text{ and } \lim_{l \to +\infty} b_{k_l+1} b_{k_l+2} \ldots b_{k_l+n_l} \frac{1}{q_{k_l+n_l}}=+\infty.$$
Let $I_{k_l-1}=\left[0,\frac{1}{b_0 b_1 \cdots b_{k_l-1}}\right]$ and $q_{k_l+n_l} I_{k_l+n_l}=\left[0,\frac{q_{k_l+n_l}}{b_0 b_1 \cdots b_{k_l+n_l}}\right]$. Then for sufficiently large $l$, we have that $\frac{\log q_{k_l} \cdots q_{k_l+n_l} \frac{1}{q_{k_l+n_l}}}{\log b_{k_l} \cdots b_{k_l+n_l} \frac{1}{q_{k_l+n_l}}}<t$ and that
$$\frac{\mu(I_{k_l-1})}{\mu(q_{k_l+n_l} I_{k_l+n_l})}\bigg/\left(\frac{|I_{k_l-1}|}{|q_{k_l+n} I_{k_l+n_l}|}\right)^t=\left(b_{k_l} \cdots b_{k_l+n_l} \frac{1}{q_{k_l+n_l}}\right)^{\frac{\log q_{k_l} \cdots q_{k_l+n_l} \frac{1}{q_{k_l+n_l}}}{\log b_{k_l} \cdots b_{k_l+n_l} \frac{1}{q_{k_l+n_l}}}-t}\to 0.$$
Hence, by definition of lower dimension, we have $\dim_L\mu\geq t$.

On the other hand, we prove $\dim_L\mu \geq s$. It is clear that $\dim_L\mu\geq 0$. Without loss of generality, we only need to consider the case that $s>0$. Consider any $r\leq R$ and $x\in \textrm{supp}(\mu)$. Choose 
integers $k, l, u \in \N$ and $n \ge 0$ such that 
$$\frac{l}{b_0 b_1 \cdots b_k} < R \leq \frac{ l+1}{b_0 b_1 \cdots b_{k}} \leq \frac{2 l}{b_0 b_1 \cdots b_{k}}, \quad 1 \le l < b_k$$
and
$$\frac{u}{b_0 b_1 \cdots b_{k+n}} < u \leq \frac{ u+1}{b_0 b_1 \cdots b_{k+n}} \leq \frac{2 u}{b_0 b_1 \cdots b_{k+n}}, \quad 1 \le u < b_{k+n}.$$ 
Consequently, 
$$\frac{\min\{l,q_k\}}{q_0 q_1 \cdots q_k}\leq \mu(B(x,R))\leq \frac{2 \min\{l+1,q_k\}}{q_0 q_1 \cdots q_k} \le \frac{4 \min\{l,q_k\}}{q_0 q_1 \cdots q_k}$$
and
$$\frac{\min\{u,q_{k+n}\}}{q_0 q_1 \cdots q_{k+n}} \leq \mu(B(x,r))\leq \frac{2 \min\{ u+1,q_{k+n}\}}{q_0 q_1 \cdots q_{k+n}} \leq \frac{4 \min\{ u,q_{k+n}\}}{q_0 q_1 \cdots q_{k+n}}.$$
This yields for any $0\leq t<s$ ($\le 1$),
\begin{equation}\label{eqleqa1}
\begin{aligned}
	& \frac{\mu(B(x,R))}{\mu(B(x,r))}\bigg/\bigg(\frac{R}{r}\bigg)^t \geq \frac{\frac{\min\{l,q_k\}}{4 \min\{u,q_{k+n}\}} \cdot q_{k+1} q_{k+2} \cdots q_{k+n}}{(\frac{2 l}{u} \cdot b_{k+1} b_{k+2} \cdots b_{k+n})^t} \\
	\ge& \begin{cases}
	\frac{1}{8} \cdot \frac{q_{k+1} q_{k+2} \ldots q_{k+n} \frac{1}{q_{k+n}}}{\left(b_{k+1} b_{k+2} \ldots b_{k+n} \frac{1}{q_{k+n}}\right)^t} & \text{if } 1 \le l < q_k; \\
	\frac{1}{8} \cdot \frac{q_{k} q_{k+1} \ldots q_{k+n} \frac{1}{q_{k+n}}}{\left(b_{k} b_{k+1} \ldots b_{k+n} \frac{1}{q_{k+n}}\right)^t} & \text{if } q_k \le l < b_k.
	\end{cases}
\end{aligned}
\end{equation}
For the above $t$, there exists $N=N(t)$ such that
$$ \frac{\log q_{k} q_{k+1} \cdots q_{k+n} \frac{1}{q_{k+n}}}{\log b_{k} b_{k+1} \cdots b_{k+n} \frac{1}{q_{k+n}}} > t \text{ whenever } {b_{k} b_{k+1} \cdots b_{k+n} \frac{1}{q_{k+n}}} > N.$$
Hence, equation \eqref{eqleqa1} is no less than $\min\{8^{-1},8^{-1}\frac{1}{N^t}\}\geq \frac{1}{8N}$. This proves $\dim_L\mu\geq t$, which implies $\dim_L\mu\geq s$ since $t > s$ is arbitrary.

Now we show that 
$$\lim_{N\to+\infty} \inf_{n \ge N, k\geq1}\frac{\log q_k q_{k+1}\cdots q_{k+n}}{\log b_k b_{k+1}\cdots b_{k+n}} = \lim_{N \to +\infty} \inf_{b_{k+1}\cdots b_{k+n}\frac{1}{q_{k+n}} \ge N}\frac{\log q_{k+1}q_{k+2}\cdots q_{k+n}\frac{1}{q_{k+n}}}{\log b_{k+1}b_{k+2}\cdots b_{k+n}\frac{1}{q_{k+n}}},$$
when $\{b_n\}_{n=1}^{\infty}$ is bounded above by $M$. Indeed, it follows from the fact
\begin{align*}
    & \{(k,n): n \ge N, k \ge 1\} \\
    \subseteq& \{(k,n): b_{k} b_{k+1}\cdots b_{k+n} \frac{1}{q_{k+n}} \ge 2^N, k \ge 1\} \\
    \subseteq& \{(k,n): n \ge N \log_M 2 - 1, k \ge 1\}.
\end{align*}

Now we assume $\{b_n\}$ is unbounded and that $\liminf_{n\to+\infty} \frac{q_n}{b_n^{1-\epsilon}}=0$ for any $\epsilon\in[0,1)$. Then, there exists a subsequence $\{b_{n_k}\}_{k=1}^{\infty}$ of $\{b_{n}\}_{n=1}^{\infty}$ such that $b_{n_k}\to +\infty$ as $k\to+\infty$. Let $I_n=\left[0,\frac{1}{b_0 b_1 \cdots b_n}\right]$ and $b_nI_n=\left[0,\frac{1}{b_0 b_1 \cdots b_{n-1}}\right]$, where $n\geq 1$. Then
$$\frac{\mu(b_{n_k}I_{n_k})}{\mu(I_{n_k})}=q_{n_k}\quad \text{while}~\left(\frac{|b_{n_k}I_{n_k}|}{|I_{n_k}|}\right)^t=(b_{n_k})^t$$
for any $t\geq0$. Therefore, when $0< t< 1$, since $\lim_{n\to\infty}\frac{q_n}{b_n^t}=0$, $$\frac{\mu(q_{n_k}I_{n_k})}{\mu(I_{n_k})}\bigg/\left(\frac{|q_{n_k}I_{n_k}|}{|I_{n_k}|}\right)^t=\frac{q_{n_k}}{(b_{n_k})^t} \to0~\text{as $k\to +\infty$}.$$
Hence, we have $\dim_L \mu=0$. 
\eproof

\section{Proofs of Theorems \ref{the-main} and \ref{thm:arbitrary_dimension}}
In this section, we shall give the proofs of Theorem \ref{the-main} and Theorem \ref{thm:arbitrary_dimension}. We first prove Theorem \ref{the-main}. Exploiting the relation \eqref{eqdim}, we only need to construct a measure $\mu$ whose lower dimension is the same as the Assouad dimension. For any $\gamma\in(0,1]$, take natural numbers $\alpha_0,\alpha_1,\beta \ge 2$ such that both $\alpha_0$ and $\alpha_1$ divide $\beta$ and that
$$\frac{\log \alpha_0}{\log \beta}<\frac{\log \alpha_0 + \log \alpha_1}{2 \log \beta}<\gamma\le\frac{\log \alpha_1}{\log \beta}. $$
Let $\Sigma=\{0,1\}$. For $n\ge 1$, write
\[
\Sigma^n=\{u=u_1 u_2 \cdots u_n: u_i \in \Sigma, 1\le i\le n\}
\]
and $\Sigma^*=\cup_{n=1}^\infty \Sigma^n$. We also write $\Sigma^\infty=\Sigma\times\Sigma\times\cdots$. Let $u=u_1u_2\cdots\in\Sigma^\ast\cup\Sigma^\infty$, define  $u_{[n,m]}=u_nu_{n+1}\cdots u_{m}$ for $1\leq n<m$.  We define a map $\Phi: \Sigma^* \to \Big[\frac{\log \alpha_0}{\log \beta},\frac{\log \alpha_1}{\log \beta}\Big]$ by assigning $u=u_1u_2\cdots u_n \in \Sigma^n$ the value
$$\Phi(u)=\frac{\#\{1\leq i \leq n:u_i=0\}\cdot\log\alpha_0+\#\{1\leq i \leq n:u_i=1\}\cdot\log\alpha_1}{n \log \beta}.$$
In the following, we will construct a sequence $x=x_1x_2\cdots\in \Sigma^\infty$ that satisfies 
$$
\liminf_{r \to \infty} \inf_{k \ge 1} \Phi(x_{[k,k+r]})=\limsup_{r \to \infty} \sup_{k \ge 1} \Phi(x_{[k,k+r]})=\gamma. 
$$
Hence, by letting $b_n=\beta$ and $q_n=\alpha_{x_n}$, the associated Moran measure $\mu$ has the property $\dim_L \mu = \dim_A \mu$. Moreover, since $\alpha_0$ and $\alpha_1$ divide $\beta$, $\mu$ is a spectral measure as is proved in \cite{AH14}.

Define a sequence $\{(u(k),v(k))\}_{k=1}^{\infty}$ of $(\Sigma^\ast)^2$ as follows. Let $u(0)=0$, $v(0)=1$, and $n_0=1$. Suppose $(u(k-1),v(k-1))$ is defined. If $k$ is odd, define 
$$
n_{k}=\max \{n \ge 1: \Phi(v(k-1)^{n-1} u(k-1)) \le \gamma\},
$$
and
\begin{equation} \label{sequence_construction}
    (u(k),v(k))=\begin{cases}
    (v(k-1)^2, v(k-1)^2) & \text{if } \Phi(v(k-1))=\gamma; \\
    (v(k-1)^{n_{k}-1} u(k-1), v(k-1)^{n_{k}}) & \text{otherwise.}\end{cases}
\end{equation}
If $k$ is even, define 
$$
n_{k}=\max \{n \ge 1: \Phi(u(k-1)^{n-1} v(k-1)) \ge \gamma\},
$$
and
\[ \label{sequence_construction*}
    (u(k),v(k))=\begin{cases}
    (u(k-1)^2, u(k-1)^2) & \text{if } \Phi(u(k-1))=\gamma; \\
    (u(k-1)^{n_{k}}, u(k-1)^{n_{k}-1} v(k-1)) & \text{otherwise.}\end{cases} \tag{\ref{sequence_construction}*}
\]
Denote by $N_k=\abs{u(k)}=\abs{v(k)}$ the length of the $k$-th word in the sequence. This choice of $(u(k),v(k))$ has the following properties.
\begin{lem}\label{lim1}
	Let $\{(u{(k)},v{(k)})\}_{k=1}^{\infty}$ be defined as above. Then, the following is true for all $k \in \Z_+$.
	\begin{enumerate}[label=(\roman*), font=\normalfont]
	    \item $0 \le \gamma - \Phi(u(k)) \le \tfrac{\Phi(v(k))-\gamma}{n_k}$ if $k$ is odd, and $0 \le \Phi(v(k))-\gamma \le \tfrac{\gamma - \Phi(u(k))}{n_k}$ if $k$ is even. In addition, $2 \le n_{k+1} \le \infty$.
	    \item $u(k),v(k) \in \{u(j),v(j)\}^{N_{k}/N_j}$ for all $0 \le j \le k$. Moreover, 
	    \[
	   u(k+1)_{[1,N_{k}]}=v(k+1)_{[1,N_{k}]}=\begin{cases} v(k), & \text{if } k+1 \text{ is odd;} \\
	    u(k), & \text{if } k+1 \text{ is even.}
	    \end{cases}
	    \]
	    \item $\lim_{k \to \infty} \Phi(u(k))=\lim_{k \to \infty} \Phi(v(k))=\gamma$. 
	\end{enumerate}
\end{lem}
\proof
(i) We prove the claim by induction on $k$. The case $k=0$ follows from the choice of $\alpha_0, \alpha_1$, and $\beta$. For induction step, assume the claim holds for $k-1$. If $k$ is odd, then by the induction hypothesis we have 
\begin{equation} \label{eq:induction_step_1}
    0 \le \Phi(v(k-1)) - \gamma \le \tfrac{\gamma - \Phi(u(k-1))}{n_{k-1}} \quad\text{ and }\quad 2 \le n_{k}  \le \infty.
\end{equation}
It is noteworthy that the proof is already finished at this point if $\Phi(v(k-1))=\gamma$, since that $\Phi(u(k))=\Phi(v(k))=\gamma$ and that $n_{k+1}=\infty$ follow immediately from \eqref{sequence_construction} in this case. Now, it remains to show that the claim also holds when $\Phi(v(k-1)) > \gamma$. Under the circumstances, $n_k$ is finite and the maximality of $n_k$ suggests that $\Phi(v(k-1)^{n_k} u(k-1)) > \gamma$ or, equivalently, $n_k > \frac{\gamma - \Phi(u(k-1))}{\Phi(v(k-1)) - \gamma}$. Therefore,
\begin{align*}
    0 & \le \gamma-\Phi(u(k)) = \gamma - \Phi(v(k-1)^{n_{k}-1} u(k-1)) \\
    & = \left.\tfrac{1}{n_{k}} \middle(\gamma - \Phi(u(k-1))\right) - \left.\tfrac{n_{k} - 1}{n_{k}} \middle(\Phi(v(k-1))-\gamma\right) \\
    & < \left.\tfrac{1}{n_k}\middle(\Phi(v(k-1))-\gamma\right) = \left.\tfrac{1}{n_k}\middle(\Phi(v(k))-\gamma\right),
\end{align*}
from which $2 \le n_{k+1} \le \infty$ follows as a consequence. This proves the case when $k$ is odd, and the argument for even $k$ is similar. Thus, the claim is valid for all $k \in \Z_+$ by induction.

(ii)  It follows immediately from \eqref{sequence_construction} and \eqref{sequence_construction*}.

(iii) Since $\Phi(v(k))=\Phi(v(k-1))$ if $k$ is odd and $\Phi(u(k))=\Phi(u(k-1))$ if $k$ is even, we deduce from \eqref{sequence_construction}, \eqref{sequence_construction*} and (i) that 
\[
0 \le \max\{\gamma-\Phi(u(k)),\Phi(v(k))-\gamma\} \le 2^{-k} (\Phi(v(0))-\gamma).
\]
The limit thus exists.
\eproof

Since Lemma \ref{lim1}(i) along with \eqref{sequence_construction} and \eqref{sequence_construction*} implies that $\lim_{k \to \infty} N_k = \infty$, Lemma \ref{lim1}(ii) naturally yields a sequence of symbols $x \in \Sigma^{\infty}$ satisfying $x_{[1,N_{2k}]}=v(2k)$ for all $k$. The following lemma demonstrates that this sequence $x$ exhibits the desired asymptotic behavior.
\begin{lem}\label{lim2}
    Let $\{(u(k),v(k))\}_{k=1}^{\infty}$ and $x$ be defined as above.  Then 
	$$\limsup_{r\to+\infty}\sup_{n \ge 1}\Phi(x_{[n,n+r]})=\liminf_{r\to+\infty}\inf_{n \ge 1}\Phi(x_{[n,n+r]})=\gamma.$$
\end{lem}
\proof
Let $k$ be fixed. For any $\epsilon > 0$, take $L \in \N$ such that $2 N_k < L \epsilon$. As a consequence of the choice of $L$, for all $n \in \N$ and all $r \ge L$, there exists a subinterval $[p N_k + 1, q N_k]$ of $[n, n+r]$ such that $r+1-(q-p) N_k \le 2 N_k < L \epsilon$ and that $p,q\in\Z_+$. Thus, it follows from Lemma \ref{lim1}(ii) that
\begin{align*}
	\Phi(x_{[n,n+r]}) &=\frac{(q-p) N_k}{r+1} \Phi(x_{[p N_k + 1, q N_k]})+\frac{r+1-(q-p) N_k}{r+1} \Phi(x_{[n,p N_k]} x_{[q N_k+1,n+r]})\\
	& \in ((1-\epsilon) \Phi(u(k)),(1-\epsilon) \Phi(v(k)) + \epsilon).
\end{align*}
Since $\epsilon$ is arbitrary, 
$$
\Phi(u(k)) \le \liminf_{r\to+\infty}\inf_{n \ge 1}\Phi(x_{[n,n+r]}) \le \limsup_{r\to+\infty}\sup_{n \ge 1}\Phi(x_{[n,n+r]}) \le \Phi(v(k))
$$
for the given $k$. Applying Lemma \ref{lim1}(iii) and letting $k \to \infty$, we deduce 
$$
\liminf_{r\to+\infty}\inf_{n \ge 1}\Phi(x_{[n,n+r]})=\limsup_{r\to+\infty}\sup_{n \ge 1}\Phi(x_{[n,n+r]})=\gamma,
$$
and the proof is completed.
\eproof

{\noindent{\bf Proof of Theorem \ref{the-main}.}}
The previous discussions prove the theorem for $t\in(0,1]$. For the case $t=0$, we only need to take some positive integers $2 \le \alpha < \beta$ such that $\alpha$ divides $\beta$ and define $q_n=\alpha, b_n=\beta^n$ for $n \in \mathbb{N}$. It is not hard to verify that $\dim_A\mu=0$ by Proposition \ref{ass}, which completes the proof.
\eproof

We now turn to the proof of Theorem \ref{thm:arbitrary_dimension}. Before we carry out the proof, we shall introduce a few notations. Given any $0\leq t_L\leq t_H\leq t_P\leq t_A\leq 1$, we may take positive integers $\alpha_0, \alpha_1, \beta \ge 2$ such that both $\alpha_0$ and $\alpha_1$ divide $\beta$ and that
\[
\frac{\log \alpha_0}{\log \beta} \le \min \left(\{t_L,t_H,t_P,t_A\} \cap (0,1]\right) \le \max \{t_L,t_H,t_P,t_A\} \le \frac{\log \alpha_1}{\log \beta}.
\]
For this choice of $\alpha_0, \alpha_1$, and $\beta$, one may define measures $\mu_\star$ determined by sequences $\{b_{\star,k}\}_{k=1}^\infty$ and $\{q_{\star,k}\}_{k=1}^\infty$ such that $\dim_\star \mu_\star = t_\star$ for $\star=L,H,P,A$ and that
\[
\log b_{\star,k}=\begin{cases}
\log \beta, & \text{if } t_\star > 0; \\
k \log \beta, & \text{otherwise},
\end{cases} \text{ and } \log q_{\star,k}\begin{cases}
\in \{\log\alpha_0,\log\alpha_1\}, & \text{if } t_\star > 0; \\
=\alpha_0, & \text{otherwise}.
\end{cases}
\]
Now we are ready to construct the desired measure $\mu$ associated with sequences $\mathrm{b}_k$ and $\mathrm{q}_k$. Take $M_{\star,k}$ to be sequences of natural numbers satisfying the following conditions:
\begin{enumerate}
    \item $M_{A,k} = M_{L,k} = \lfloor k^{1/2} \rfloor$, where  $\lfloor x\rfloor$ is the greatest integer which is not larger than $x$, and
    \item $M_{P,k} \ge k M_{H,k}^2$ and $M_{H,k+1} \ge (k+1) M_{P,k}^2$.
\end{enumerate}
Define $\mathrm{b}_k$ and $\mathrm{q}_k$ as
\begin{equation*}
\begin{aligned}
\{\mathrm{b}_k\}_{k=1}^\infty=\{&b_{H,1},\ldots,b_{H,M_{H,1}},    b_{P,1},\ldots,b_{P,M_{P,1}},    b_{L,1},\ldots,b_{L,M_{L,1}},    b_{A,1},\ldots,b_{A,M_{A,1}}, \\
&b_{H,1},\ldots,b_{H,M_{H,2}},    b_{P,1},\ldots,b_{P,M_{P,2}},    b_{L,1},\ldots,b_{L,M_{L,2}},    b_{A,1},\ldots,b_{A,M_{A,2}} ,\ldots\},
\end{aligned}
\end{equation*}
and
\begin{equation*}
\begin{aligned}
\{\mathrm{q}_k\}_{k=1}^\infty=\{&q_{H,1},\ldots,q_{H,M_{H,1}},    q_{P,1},\ldots,q_{P,M_{P,1}},    q_{L,1},\ldots,q_{L,M_{L,1}},    q_{A,1},\ldots,q_{A,M_{A,1}}, \\
&q_{H,1},\ldots,q_{H,M_{H,2}},    q_{P,1},\ldots,q_{P,M_{P,2}},    q_{L,1},\ldots,q_{L,M_{L,2}},    q_{A,1},\ldots,q_{A,M_{A,2}} ,\ldots\}.
\end{aligned}
\end{equation*}
For simplicity, we denote the index of $b_{\star,M_{\star,k}}$ in $\{b_k\}_{k=1}^{\infty}$ as $S_{\star,k}$, and the interval $[S_{\star,k}-M_{\star,k}+1,S_{\star,k}]$ as $I_{\star,k}$. In other words, $S_{H,1}=M_{H,1}$, $S_{P,1}=M_{H,1}+M_{P,1}$, and so on. Here and in the following, we write interval $[n,m]$ as $\{n,n+1,\ldots, m\}$ for $n<m$. 
 
 For any interval $J \subseteq \mathbb{N}$, define $\rho(J)=\Pi_{i \in J} \mathrm{q}_i$ and $\lambda(J)=\Pi_{i \in J} \mathrm{b}_i$. Note that since $M_{\star,k} \le \frac{\log\lambda(I_{\star,k})}{\log\beta} \le M_{\star,k}^2$, it follows from (1) and (2) that
\begin{enumerate}
    \item[(1')] $\frac{\log\lambda(I_{A,k})}{\log\beta} = \frac{\log\lambda(I_{L,k})}{\log\beta} \le k$, and
    \item[(2')] $\frac{\log\lambda(I_{P,k})}{\log\beta} \ge \frac{\log\lambda(I_{H,k})}{\log\beta}$ and $\frac{\log\lambda(I_{H,k+1})}{\log\beta} \ge (k+1) \frac{\log\lambda(I_{P,k})}{\log\beta}$,
\end{enumerate}
which yield the following lemma. We omit the proof since it is straightforward.
\begin{lem} \label{lem:estimate_density}
    The following statements hold.
    \begin{enumerate}[label=(\roman*), font=\normalfont]
        \item $\lim_{k \to \infty} \frac{\log\lambda(I_{\star,k})}{\log\lambda([1,S_{\star,k}])}$ is $1$ if $\star=H,P$ and is $0$ if $\star=L,A$.
        \item For all $R \in \mathbb{N}$, $\liminf_{\lambda(J) \to \infty} \tfrac{\sum_{\star,k: \lambda(I_{\ast,k} \cap J) \ge R}\log\lambda(I_{\ast,k} \cap J)}{\log\lambda(J)}=1$.
    \end{enumerate}
\end{lem}
{\noindent{\bf Proof of Theorem \ref{thm:arbitrary_dimension}.}}
Since $\mathrm{q}_k$ is bounded, the dimension formulae are reduced to 
\begin{align*}
    \dim_L\mu=\liminf_{\lambda(J) \to \infty} \frac{\log \rho(J)}{\log \lambda(J)}, \quad \dim_H\mu=\liminf_{n \to \infty} \frac{\log \rho([1,n])}{\log \lambda([1,n])}, \\
    \dim_P\mu=\limsup_{n \to \infty} \frac{\log \rho([1,n])}{\log \lambda([1,n])}, \quad \dim_A\mu=\limsup_{\lambda(J) \to \infty} \frac{\log \rho(J)}{\log \lambda(J)},
\end{align*}
where the formulae of Hausdorff dimension and packing dimension of $\mu$ are given in \cite{FWW}.

We first show that $\dim_A \mu=t_A$. For $\dim_A \mu \ge t_A$, we note that 
\[
\dim_A \mu = \limsup_{\lambda(J) \to \infty} \frac{\log \rho(J)}{\log \lambda(J)} \ge \lim_{k \to \infty} \frac{\log \rho(I_{A,k})}{\log \lambda(I_{A,k})} = t_A. 
\]
On the other hand, the construction of $\mu_\star$ and Theorem \ref{the-main} imply that for every $\epsilon > 0$, there exists $R \in \mathbb{N}$ such that $\frac{\log q_{\star,n}\cdots q_{\star,n+r}}{\log b_{\star,n}\cdots b_{\star,n+r}} < t_A+\epsilon$ whenever $\star=L,H,P,A$ and $b_{\star,n}\cdots b_{\star,n+r} \ge R$. It then follows from Lemma \ref{lem:estimate_density}(ii) that
\begin{align*}
    & \quad \dim_A \mu = \limsup_{\lambda(J) \to \infty} \frac{\log \rho(J)}{\log \lambda(J)} \\
    &= \limsup_{\lambda(J) \to \infty}\Bigg( \sum_{\star,k:\lambda(I_{\star,k} \cap J) \ge R} \frac{\log\lambda(I_{\star,k} \cap J)}{\log\lambda(J)} \cdot \frac{\log \rho(I_{\star,k} \cap J)}{\log \lambda(I_{\star,k} \cap J)} \\
    &\hspace{7em}+\sum_{\star,k:\lambda(I_{\star,k} \cap J) < R} \frac{\log\lambda(I_{\star,k} \cap J)}{\log\lambda(J)} \cdot \frac{\log \rho(I_{\star,k} \cap J)}{\log \lambda(I_{\star,k} \cap J)}\Bigg) \\
    &\le t_A+\epsilon.
\end{align*}
Since $\epsilon$ is arbitrary, $\dim_A \mu \le t_A$. It is noted that the proof for lower dimension $\dim_L \mu = t_L$ is similar and is therefore omitted.

Now we show that $\dim_H \mu = t_H$. On one hand, it follows from Lemma \ref{lem:estimate_density}(i) that
\[
\dim_H \mu = \liminf_{\lambda(J) \to \infty} \frac{\log \rho(J)}{\log \lambda(J)} \le \lim_{k \to \infty} \frac{\log \rho([1,S_{H,k}])}{\log \lambda([1,S_{H,k}])} = t_H. 
\]
On the other, the construction of $\mu_\star$ and Theorem \ref{the-main} imply that for every $\epsilon > 0$, there exists $R \in \mathbb{N}$ such that $\frac{\log q_{\star,n}\cdots q_{\star,n+r}}{\log b_{\star,n}\cdots b_{\star,n+r}} > t_H-\epsilon$ whenever $\star=H,P$ and $b_{\star,n}\cdots b_{\star,n+r} \ge R$. It then follows from Lemma \ref{lem:estimate_density}(ii) that
\begin{align*}
    & \quad \dim_H \mu = \liminf_{n \to \infty} \frac{\log \rho([1,n])}{\log \lambda([1,n])} \\
    &= \liminf_{n \to \infty}\Bigg( \sum_{\substack{\star=H,P \\ \lambda(I_{\star,k} \cap [1,n]) \ge R}} \frac{\log\lambda(I_{\star,k} \cap [1,n])}{\log\lambda([1,n])} \cdot \frac{\log \rho(I_{\star,k} \cap [1,n])}{\log \lambda(I_{\star,k} \cap [1,n])} \\
    &\hspace{7em}+\sum_{\substack{\star,k:\lambda(I_{\star,k} \cap [1,n]) < R }} \frac{\log\lambda(I_{\star,k} \cap [1,n])}{\log\lambda([1,n])} \cdot \frac{\log \rho(I_{\star,k} \cap [1,n])}{\log \lambda(I_{\star,k} \cap [1,n])} \\
    &\hspace{7em}+\sum_{\substack{\star=L,A \\ \lambda(I_{\star,k} \cap [1,n]) \geq R }} \frac{\log\lambda(I_{\star,k} \cap [1,n])}{\log\lambda([1,n])} \cdot \frac{\log \rho(I_{\star,k} \cap [1,n])}{\log \lambda(I_{\star,k} \cap [1,n])} \Bigg)\\
    &\ge t_H-\epsilon.
\end{align*}
Since $\epsilon$ is arbitrary, $\dim_H \mu \ge t_H$. It is noted that the proof for packing dimension $\dim_P \mu = t_P$ is similar and is therefore omitted. The proof is then finished.
\eproof

	\subsection*{Acknowledgements}
This work is supported by the Academy of Finland, project grant No. 318217. We would like to thank Meng Wu for his guidance, comments, encouragement and drawing our attention to the packing dimensions of measures when preparing this paper. 

\bibliographystyle{amsplain}
\bibliography{reference}
\end{document}